\documentclass{amsart}
\usepackage{amsmath} 
\usepackage{amssymb}
\usepackage{mathrsfs}

\newtheorem{theorem}{Theorem}[section] 

\newtheorem{lemma}[theorem]{Lemma}

\theoremstyle{definition}

\theoremstyle{remark}

\newcommand{\cI}{{\mathscr I}}
\newcommand{\cK}{{\mathscr K}}

\newcommand{\mbR}{{\mathbb R}}

\newcount\skewfactor
\def\mathunderaccent#1#2 {\let\theaccent#1\skewfactor#2
\mathpalette\putaccentunder}
\def\putaccentunder#1#2{\oalign{$#1#2$\crcr\hidewidth
\vbox to.2ex{\hbox{$#1\skew\skewfactor\theaccent{}$}\vss}\hidewidth}}

\begin{document}

\title {Some corrections to an old paper}

\author {\bf E. Bombieri, J.B. Friedlander \and H. Iwaniec}

\smallskip


\thanks{JF supported in part by NSERC grant A5123}




\maketitle
 {\bf Abstract:} 
We give some corrections of our paper ``Primes in arithmetic progressions to 
large moduli'' [BFI]. The corrections do not affect the statements of any of 
the theorems in the paper. The contents of our two sequel papers [BFI2, BFI3] 
also remain unchanged.  

\section{\bf Acknowledgements} 

We are grateful to Kevin Broughan for quite recently drawing our 
attention to a slip in one 
of the theorems, Theorem 12, of the paper [DI]. We are also pleased to thank 
J-M. Deshouillers for encouraging us to write this note. We discuss the changes 
to [BFI] necessitated by this problem in Section 2. 

\medskip 

We are also grateful to Zaizhao Meng for pointing out to us (actually, quite 
some time ago) a difficulty in a separation of variables argument we employed 
in two places in the paper. We correct this in Section 3. 

\section{\bf Sums of Kloosterman fractions}

The first change concerns a result, to which reference was frequently
made in our work [BFI] and 
which is crucial to that work, recorded therein as Lemma 1 and occurring on 
page 210. It gives a bound for certain multiple sums of Kloosterman fractions. 
We now 
state it as follows:

\begin{lemma}
Let $g_0(\xi,\eta)$ be a smooth function with compact support in 
$\mbR^+ \times \mbR^+ $. Let $C, D, N, R, S > 0$ and $g(c,d) = g_0(c/C, d/D)$. 
For any complex numbers $B_{nrs}$ denote 

$$
\cK(C,D,N,R,S) = \sum_{r\sim R}\sum_{s\sim S}\sum_{0<n\leqslant N}B_{nrs}
\mathop{\sum_c\sum_d}_{(rd,sc)=1} 
g(c,d) e \left(n\frac{\overline{rd}}{sc}\right)\ . 
$$
Then, for any $\varepsilon > 0$ we have 
$$
\cK(C,D,N,R,S) \ll (CDNRS)^{\varepsilon} \cI (C,D,N,R,S)\, \Vert B \Vert\ ,
$$
where $\Vert B \Vert =\Vert B \Vert_2$ is the $\ell_2$-norm and 
$\cI^2(C,D,N,R,S)$ is the quantity 
$$
 CS(RS+N)(C+DR) + C^2DS\sqrt{(RS+N)R} +D^2NR\ ,
$$
the constant implied in $\ll$ depending at most on $\varepsilon$ and 
$g(\xi,\eta)$. 

\end{lemma}

This is somewhat weaker than the version quoted in [DI] and [BFI] 
where, in the final term, the quantity $D^2NR$ was 
stated as $D^2NRS^{-1}$. The above corrected version is already sufficient
for our applications. This modification of the final term is the follow-up
of the correction of the bound (9.11) of [DI], wherein the quantity
$D(NR/S)^{1/2}$ needs to be replaced by $D(NR)^{1/2}$. No other changes in [DI]
need to be made as a consequence of this replacement.

In most of our uses of this lemma we have $S=1$ so things remain as before. 
This is not however the case in Section 10 where there is needed a small 
change in 
one line in the proof (but not in the statement) of Lemma 8. This occurs in the 
first display following (10.8) on page 233, the final one in the proof of 
Lemma 8. In the final term of that display $Q^2HN^3$ needs to be replaced by 
the larger term $Q^2HN^5$. Nothing further needs to be changed however because 
even this larger version of the last term is dominated by at least one of the 
preceding terms, for example by the term $Q^2HN^7$.

\section {\bf A separation of variables argument}

The argument in question occurs in Section 9 of the paper on pages 226 and 227 
and the formula (9.12) there is incorrect. The same argument 
is repeated in very slightly different form in Section 11 on page 234. 
The goal of each discussion, namely (9.13), (11.3), remains valid, 
with a value of $K$ only slightly larger than that 
used in the paper (see (9.6)). The reason why this is so is that $\delta$ is 
a divisor of $a$, which is assumed to be fixed throughout the paper, and 
$q_0$ is small. In fact $q_0\le {\mathcal L}^{A+B}$, so it plays only a 
negligible role in the estimate.  What we have to do is to work directly with 
the modulus $\delta q_0 k$, rather than separately with the moduli 
$\delta q_0$ and $k$ as in the paper. The corrections to be made are as follows.

\bigskip

\noindent{\it p. 226, from (9.6) up to (9.7), replace with:} 
$$ 1\leq |k| \leq K_0\eqno (9.6)$$
where now $K_0=N/q_0R$.

We wish to separate the variables $h$, $n_2$ from the remaining ones. We detect
the conditions (9.4) by means of multiplicative characters 
$\chi\,
({\rm mod}\,\,
\delta q_0 k)$, that is we appeal to the following orthogonality relation
$${1\over\varphi(\delta q_0k)}
\sum_{\chi\,({\rm mod}\,\,
\delta q_0 k)}\overline\chi(n_1)\chi(n_2)=
\begin{cases} 
1 \quad \quad {\rm if}\,\,\,  n_1\equiv n_2\, ({\rm mod}\,\,
\delta q_0 k),\  (n_1n_2,\delta q_0 k)=1\cr
0\quad \quad\quad\quad{\rm otherwise.}\quad \quad \quad \quad\quad\quad \hfill
\rm(9.7)
\end{cases}
$$
\medskip
\noindent{\it p. 227, formula (9.12), replace with:}
\begin{align*}
{\mathcal R}_1  \ll &\, Y\,(\log 2N)\sum_{\delta | a}\sum_{q_0\le 
Q_0}\sum_{1\leq k\leq K_0} {1\over\varphi(\delta q_0 
k)}\sum_{\chi\,({\rm mod}\,\,\delta q_0 k)}\,
\sum\sum_{\kern-16pt (q_1,q_2)=1}|\gamma_{q_0q_1}\gamma_{q_0q_2}|\cr
\quad & \times\sum_{(n_1,q_1)=1}|\beta_{n_1}|\sum_{1\leq |h|\leq H}
\left|\sum_{(n_2,n_1q_2)=1}\beta(h,n_2) \chi(n_2) 
\,e\!\left(ahk{\overline{n_2q_1}\over n_1q_2}\right)\right| .
\quad \quad \,\, \rm{(9.12)}
\end{align*}

\medskip
\noindent{\it p. 228, line 2, replace with:}

\noindent with some coefficients $\beta(h,n)$ such that $|\beta(h,n)|\leq 
|\beta_n|$.  

\medskip 
\noindent{\it p. 234, last display before (11.3), replace with:}
\begin{align*}
{\mathcal R}_1  \ll \, x^{\varepsilon}\,MR^{-1}&\sum_{\delta | a}\sum_{1\leq k\leq K_0} {1\over\varphi(\delta k)}\sum_{\chi\,({\rm mod}\,\,\delta k)}\,
\sum_{n_1}\sum_{l_1}\sum_{l_2}|\beta_{n_1}\lambda_{l_1}\lambda_{l_2}|\cr
\quad \quad\quad \quad \quad\quad & \times\left|\sum_{1\leq |h|\leq H}
\sum_{n_2}\beta(h,n_2) \chi(n_2) 
\,e\!\left(ahk{\overline{n_2l_2}\over n_1l_1}\right)\right| .
\end{align*}

\bigskip

In order to see that these suffice for our purpose, note that the new value of 
$K$ in (9.14), respectively (11.3), is at most $|a|Q_0$ times, respectively $|a|$ times, the old value of $K$ (which now become $K_0$), 
while $|a|Q_0\ll x^\varepsilon$ for any fixed $a$ and $\varepsilon>0$.  Since 
all estimates in the rest of Sections 9 and 11 allow for this factor 
$x^\varepsilon$, no further changes are needed.

\medskip

 In conclusion, we mention that the corrections in this section were
made in 2001 and sent in response, shortly after receipt of the
communication from Professor Meng. They were also submitted, but not 
published, at that time "since the paper 
is old and the corrections are not vital".

\bigskip
\medskip 
School of Mathematics, Institute for Advanced Study

Princeton, NJ 08540, USA

\medskip 
Department of Mathematics, University of Toronto

Toronto, Ontario M5S 2E4, Canada 

\medskip

Department of Mathematics, Rutgers University

Piscataway, NJ 08903, USA


\begin{thebibliography}{xxxx}

\bibitem[BFI]{BFI} E. Bombieri, J.B. Friedlander and H. Iwaniec, Primes in 
arithmetic progressions to large moduli, {\it Acta Math.} {\bf 156} (1986) 
203--251.

\bibitem[BFI2]{BFI2} E. Bombieri, J.B. Friedlander and H. Iwaniec, Primes in 
arithmetic progressions to large moduli II,  {\it Math. Ann.} {\bf 277} (1987) 
361--393.

\bibitem[BFI3]{BFI3} E. Bombieri, J.B. Friedlander and H. Iwaniec, Primes in 
arithmetic progressions to large moduli III, {\it J. Amer. Math. Soc.} 
{\bf 2} (1989)  215--224.

\bibitem[DI]{DI} J-M. Deshouillers and H. Iwaniec, Kloosterman sums and 
Fourier coefficients of cusp forms, {\it Invent. Math.} {\bf 70} (1982/83) 
219--288. 



\end{thebibliography}
\end{document}